\RequirePackage{ifpdf}
\ifpdf 
\documentclass[pdftex]{sigma}
\else
\documentclass{sigma}
\fi

\numberwithin{equation}{section}

\newcommand{\lltwoz}{\ell^2({\mathbb Z})}

\newcommand{\wlltwoz}{\ell_{a}^2({\mathbb Z})}

\newcommand{\D}{{\mathbb D}}

\newcommand{\R}{{\mathbb R}}
\newcommand{\C}{{\mathbb C}}

\newcommand{\Z}{{\mathbb Z}}

\renewcommand{\d}{\partial}

\newcommand{\f}{\varphi}

\begin{document}

\allowdisplaybreaks

\renewcommand{\thefootnote}{$\star$}

\renewcommand{\PaperNumber}{056}

\FirstPageHeading

\ShortArticleName{A Note on Dirac Operators on the Quantum Punctured Disk}

\ArticleName{A Note on Dirac Operators\\ on the Quantum Punctured Disk\footnote{This paper is a
contribution to the Special Issue ``Noncommutative Spaces and Fields''. The
full collection is available at
\href{http://www.emis.de/journals/SIGMA/noncommutative.html}{http://www.emis.de/journals/SIGMA/noncommutative.html}}}

\Author{Slawomir KLIMEK and Matt MCBRIDE}

\AuthorNameForHeading{S.~Klimek and M.~McBride}

\Address{Department of Mathematical Sciences,
Indiana University-Purdue University Indianapolis,\\
402 N. Blackford St., Indianapolis, IN 46202, USA}
\Email{\href{mailto:sklimek@math.iupui.edu}{sklimek@math.iupui.edu}, \href{mailto:mmcbride@math.iupui.edu}{mmcbride@math.iupui.edu}}

\ArticleDates{Received March 30, 2010, in f\/inal form July 07, 2010;  Published online July 16, 2010}

\Abstract{We study quantum analogs of the Dirac type operator $-2\overline{z}\frac{\d}{\d\overline{z}}$ on the punctured disk, subject to the Atiyah--Patodi--Singer boundary conditions.  We construct a parametrix of the quantum operator and show that it is bounded outside of the zero mode.}


\Classification{46L99; 47B25; 81R60}

\renewcommand{\thefootnote}{\arabic{footnote}}
\setcounter{footnote}{0}

\section{Introduction}

The main technical and computational part of the Atiyah, Patodi, Singer paper \cite{APS} is the initial section containing a study of a nonlocal boundary value problem for the f\/irst order dif\/ferential operators of the form $\Gamma(\frac{\d}{\d t}+B)$ on the semi-inf\/inite cylinder
$\R_+\times Y$, where $t\in\R_+$ and~$B$,~$\Gamma$ live on the boundary $Y$. The novelty of the paper was the boundary condition, now called the APS boundary condition, that involved a spectral projection of $B$. The authors explicitly compute and estimate the fundamental solutions on the cylinder. This is later used to construct a parametrix for the analogical boundary value problem on a manifold with boundary by gluing it with a contribution from the interior, see also \cite{BBW}.

The present paper aims, in a special case, to reproduce such results in the noncommutative setup of \cite{Connes}. A similar but dif\/ferent study of an example of APS boundary conditions in the context of noncommutative geometry is contained in \cite{CPR}.

This paper is a continuation of the analysis started in \cite{CKW} and \cite{KM}. The goal of those articles was to provide simple examples of Dirac type operators on noncommutative compact manifolds with boundary and then study Atiyah--Patodi--Singer type boundary conditions and the corresponding index problem.
This was done for the noncommutative disk and the noncommutative annulus and for two somewhat dif\/ferent types of operators constructed by taking commutators with weighted shifts.

In this paper we consider such non-commutative analogs of the Dirac type operator $\frac{\d}{\d t} + \frac{1}{i}\frac{\d}{\d\f}$ on the cylinder $\R_+\times S^1$, which we view as a punctured disk.  Using a weighted shift, which plays the role of the complex coordinate $z$ on the disk, we construct quantum Dirac operators, and analogs of the $L^2$ Hilbert space of functions in which they act.
We then consider the boundary condition of Atiyah, Patodi, Singer. This is done in close analogy with the commutative case. The main result of this note is that a quantum operator has an inverse which, minus the zero mode, is bounded just like in Proposition~2.5 of \cite{APS}.  In contrast with our previous papers the analysis here is more subtle because of the noncompactness of the cylinder. In particular the components of a parametrix are not compact operators and we use the Schur--Young inequality to estimate their norms. It is hoped that in the future such results will be needed to construct spectral triples and a noncommutative index theory of quantum manifolds with boundary.

The paper is organized as follows. In Section~\ref{section2} the classical APS result for the operator $-2\overline{z}\frac{\d}{\d\overline{z}}$ on the cylinder is stated and re-proved using the Schur--Young inequality. Section~\ref{section3} contains the construction of the non-commutative punctured disk and the f\/irst type of noncommutative analogs of the operator from the previous section. The operators here are similar to those of~\cite{CKW}. Also in this section a non-commutative Fourier decomposition of the Hilbert spaces and the operators is discussed.   Section~\ref{section4} contains the construction and the analysis of the Fourier components of the parametrix and the proof of the main result. Finally in Section~\ref{section5} we consider the ``balanced" versions of the quantum Dirac operators in the spirit of \cite{KM} and show how to modify the previous arguments to estimate the parametrix.

\section{Classical Dirac operator on the punctured disk}\label{section2}

In this section we revisit the analysis of Atiyah, Patodi, Singer in the simple case of semi-inf\/inite cylinder $\R_+\times S^1$, or equivalently a punctured disk. Using complex coordinates of the latter, we construct a parametrix of a version of the d-bar operator and prove norm estimates on its components by using dif\/ferent techniques than those in \cite{APS}.

Let $\D^* = \left\{z\in\C   : \; 0<|z|\le 1\right\}$ be the punctured disk.   Consider the following Dirac type operator
on $\D^*$:
\begin{gather*}
D = -2\overline{z}\frac{\d}{\d\overline{z}}.
\end{gather*}

In polar coordinates $z=re^{i\f}$ the operator $D$ has the following representation:
\begin{gather*}
D =  -r\frac{\d}{\d r}+\frac{1}{i}\frac{\d}{\d\f}=-r\frac{\d}{\d r}+B,
\end{gather*}
where $B = \frac{1}{i}\frac{\d}{\d\f}$ is the boundary operator.

We wish to study $D$, subject to the APS boundary condition, on the Hilbert space $L^2(\D^*, d\mu)$ with measure $\mu(z)$ given by the following formula:
\begin{gather}\label{measuredef}
d\mu(z) = \frac{1}{2i|z|^2}dz \wedge d\overline{z}.
\end{gather}

Let us recall the APS condition.
Let us def\/ine $P_{\ge 0}$ to be the spectral projection of~$B$ in~$L^2(S^1)$  onto the non-negative part of the spectrum of~$B$.   Equivalently, $P_{\ge 0}$ is the orthogonal projection onto $\textrm{span}\{e^{in\f}\}_{n\ge 0} $.   Then we say that $D$ satisf\/ies the APS boundary condition when its domain consists of those functions $f(z) = f(r,\f)$  on $\D^*$
which have only negative frequencies at the boundary, see \cite{APS} and \cite{BBW} for more details. More precisely:
\begin{gather}\label{ddomain1}
\textrm{dom}(D) = \left\{f\in L^2(\D^*, d\mu) \, : \   Df\in L^2(\D^*, d\mu),\    P_{\ge 0}f(1,\cdot)=0\right\}.
\end{gather}

Notice that, by the change of variable, $t=-\ln{r}$, the Dirac operator, $D$ on $L^2(\D^*, d\mu)$, is equivalent to the operator, $\frac{\d}{\d t} + \frac{1}{i}\frac{\d}{\d\f}$ on $L^2(\R_+\times S^1)$, since one has:
\begin{gather*}
d\f \wedge dt = \frac{1}{2i|z|^2}dz \wedge d\overline{z}.
\end{gather*}
This matches the APS setup.

We proceed the in the same way as in \cite{APS} by considering the spectral decomposition of the boundary operator $B$, which in our case amounts to Fourier decomposition:
\begin{gather}\label{comfour}
f(z) = \sum_{n\in\Z}f_n(r)e^{-in\f} .
\end{gather}
This yields the following decomposition of the Hilbert space $L^2(\D^*, d\mu)$:
\begin{gather}\label{hdecomp1}
L^2(\D^*, d\mu)= \bigoplus_{n\in\Z} \left(L^2\left((0,1],\frac{dr}{r}\right) \otimes \left[e^{-in\f}\right]\right)
\cong \bigoplus_{n\in\Z} L^2\left((0,1],\frac{dr}{r}\right).
\end{gather}

Now we consider the decomposition of $D$ and its inverse.  The theorem below is a special case of Proposition~2.5 of~\cite{APS} but we supply a proof that generalizes to the noncommutative setup. Def\/ine  $\overline{A}^{(n)}f(r):=-rf'(r) - nf(r)$ on the maximal domain in $ L^2((0,1],\frac{dr}{r})$, and let $\overline{A_0}^{(n)}$ be the operator $\overline{A}^{(n)}$ but with domain $\{f(r)\in \textrm{dom}(\overline{A}^{(n)}) : f(1) = 0 \}$. We have:

\begin{theorem}\label{theo1}
Let $D$ be the Dirac operator defined above on the domain \eqref{ddomain1}. With respect to the decomposition \eqref{hdecomp1} one has
\begin{gather*}
D \cong \bigoplus_{n>0} \overline{A}^{(n)} \oplus \bigoplus_{n\leq 0} \overline{A_0}^{(n)}.
\end{gather*}
Moreover, there exists an operator $Q$ such that $DQ=I=QD$, and
\begin{gather*}
Q = \bigoplus_{n\in\Z} Q^{(n)}=Q^{(0)}+\tilde{Q},
\end{gather*}
where $\tilde{Q}$ is bounded.
\end{theorem}

\begin{proof}
Staring with a function $g(z)\in L^2(\D^*,d\mu)$ we want to solve the following equation
\begin{gather*}
Df(z) = g(z)
\end{gather*}
with $f(z)$ satisfying the APS boundary condition.
The Fourier decomposition (\ref{comfour}) yields
\begin{gather*}
\sum_{n\in\Z} \left(-rf_n'(r) - nf_n(r)\right)e^{-in\f} = \sum_{n\in\Z} g_n(r)e^{-in\f}.
\end{gather*}
Therefore we must solve the dif\/ferential equation $-rf_n'(r) - nf_n(r) = g_n(r)$ where additionally $f_n(1)=0$ for $n\leq 0$.   This, and the requirement that $f$ is square integrable, assures that there is a unique solution given by the following formula:
\begin{gather*}
f_n(r)=Q^{(n)}g_n(r) = \left\{
\begin{array}{ll}
\displaystyle -\int_0^r \left(\frac{\rho}{r}\right)^n g_n(\rho)\frac{d\rho}{\rho}, \quad & n> 0,\vspace{1mm} \\
\displaystyle \int_r^1 \left(\frac{\rho}{r}\right)^n g_n(\rho)\frac{d\rho}{\rho}, \quad & n\le 0.
\end{array}\right.
\end{gather*}
This gives us the formula for the parametrix: $Q = \oplus_{n\in\Z} Q^{(n)}$.
Showing $QD=DQ=I$ is a~simple computation and is omitted.

We want to prove that $\tilde{Q}=\oplus_{n\ne 0} Q^{(n)}$ is bounded.   One has
\begin{gather*}
\|\tilde{Q}\| \le \sup_{n\in\Z\setminus\{0\}}  \big\|Q^{(n)}\big\|.
\end{gather*}
In what follows we show that the $Q^{(n)}$ are uniformly bounded, in fact of order $O\big(\frac{1}{|n|}\big)$. The main tool is the following inequality, see~\cite{HS}.

\begin{lemma}[Schur--Young inequality]\label{shuryounginq}
Let $T: L^2(Y) \longrightarrow L^2(X)$ be an integral operator:
\begin{gather*}
Tf(x) = \int K(x,y)f(y)dy.
\end{gather*}
Then one has
\begin{gather*}
\|T\|^2 \le \left(\sup_{x\in X} \int_Y |K(x,y)|dy\right)\left(\sup_{y\in Y} \int_X |K(x,y)|dx\right).
\end{gather*}
\end{lemma}

For negative $n$ one can rewrite $Q^{(n)}$ as
\begin{gather*}
Q^{(n)}g_n(r) = \int_0^1 K(r,\rho)g_n(\rho)\frac{d\rho}{\rho}
\end{gather*}
with integral kernel $K(r,\rho) = \chi(r/\rho) (r/\rho)^{|n|}$. Here the characteristic function $\chi(t)=1$ for $t\leq 1$ and is zero otherwise. Next we estimate:
\begin{gather*}
\sup_{r} \int_r^1 \frac{r^{|n|}}{\rho^{|n|+1}}d\rho = \sup_{r} \frac{1}{|n|}\big(1-r^{|n|}\big)\le \frac{1}{|n|}.
\end{gather*}
Similarly one has:
\begin{gather*}
\sup_{\rho} \int_0^\rho \frac{r^{|n|-1}}{\rho^{|n|}}dr = \sup_{\rho} \frac{1}{|n|}\cdot\frac{1}{\rho^{|n|}}\cdot\rho^{|n|} =\frac{1}{|n|}.
\end{gather*}
Thus one has by the Schur--Young inequality that $\|Q^{(n)}\|\le \frac{1}{|n|}$.   A similar computation for positive $n$ gives $\|Q^{(n)}\|\le \frac{1}{|n|}$ for all $n\neq 0$.   Hence one has that $\tilde{Q}$ is bounded.
\end{proof}

Since in this note we do not attempt to go beyond the analysis on the semi-inf\/inite cylinder, we will simply ignore $n=0$. In \cite{APS} this was not an issue as $Q^{(0)}$ is continuous when mapping into an appropriate local Sobolev space.

\section{Dirac operators on the quantum punctured disk}\label{section3}

In this section  we construct  the non-commutative punctured disk and the quantum analog of the Dirac operator of the previous section. In particular, a non-commutative Fourier decomposition of that operator is discussed. Let us also mention that a version of a quantum punctured disk was previously considered in~\cite{KL3}.

We start with def\/ining several auxiliary objects needed for our construction. Let $\lltwoz$ be the Hilbert space of square summable bilateral sequences, and let $\{e_k\}_{k\in\Z}$ be its canonical basis.   We need the following two operators: let $U$ be the shift operator given by:
\begin{gather*}
Ue_k = e_{k+1}
\end{gather*}
and let $K$ be the label operator def\/ined by the following formula:
\begin{gather*}
Ke_k = ke_k.
\end{gather*}
By the functional calculus, if $f:\Z\to\C$, then $f(K)$ is a diagonal operator and satisf\/ies the relation $f(K)e_k = f(k)e_k$.

Next assume we are given a sequence $\{w(k)\}_{k\in\Z}$ of real numbers with the following properties:
\begin{gather}\label{weightconditions}
\begin{aligned}
&1) \ w(k)<w(k+1); \\
&2) \ \lim_{k\to\infty} w(k) =: w_+ \textrm{    exists}; \\
&3) \ \lim_{k\to -\infty} w(k) = 0; \\
&4) \ \sup_{k} \frac{w(k)}{w(k-1)} < \infty.
\end{aligned}
\end{gather}
In particular we have $w(k)>0$.

The function $w:\Z\to\C$ gives a diagonal operator $w(K)$ as above.   From this we def\/ine the weighted shift operator $U_w := Uw(K)$ which plays the role of a noncommutative complex coordinate on the punctured disk.

Clearly:
\begin{gather*}
U_we_k =w(k)e_{k+1},\qquad
U^*_we_k =w(k-1)e_{k-1}.
\end{gather*}
Consider the commutator $S:= [U^*_w, U_w]$, for which one has $Se_k = (w^2(k) - w^2(k-1))e_k$.   If we let $S(k) := w^2(k) - w^2(k-1)$, then we can write $S = S(K)$.   Notice that $S$ is a trace class operator and a simple computation gives $\textrm{tr}(S) = w^2_+$.

The quantum punctured disk $C^*(U_w)$ is def\/ined to be the $C^*$-algebra generated by $U_w$.   General theory, see \cite{Conway}, gives us the following short exact sequence:
\begin{gather*}
0 \longrightarrow \mathcal{K} \longrightarrow C^*(U_w) \overset{\sigma}{\longrightarrow} C\big(S^1\big) \longrightarrow 0,
\end{gather*}
where $\mathcal{K}$ is the ideal of compact operators and $\sigma$ is the noncommutative ``restriction to the boundary'' map.

Let $b\in C^*(U_w)$ and we consider the densely def\/ined weight on $C^*(U_w)$ by
\begin{gather*}
\tau(b) = \textrm{tr}\big(S\big(U^*_wU_w\big)^{-1}b\big)
\end{gather*}
(compare with~\eqref{measuredef}). We use this weight to def\/ine the Hilbert space $\mathcal{H}$ on which the Dirac operator will live. This is done by
the GNS construction for the algebra $C^*(U_w)$ with respect to~$\tau$.   In other words $\mathcal{H}$ is obtained as a Hilbert space completion
\begin{gather*}
\mathcal{H} = \overline{(C^*(U_w), \langle\cdot, \cdot \rangle_{\tau}= \|\cdot\|^2_w)},
\end{gather*}
where $\|b\|^2_w = \tau(bb^*)$.

Now we are ready to def\/ine the operator that we wish to study, the quantum analog of the operator of the previous section.   Def\/ine $D$ by the following formula:
\begin{gather}\label{ddef}
Db = -S^{-1}U^*_w\left[b,U_w\right].
\end{gather}
Let, as before, $P_{\ge 0}$ be the orthogonal $L^2$ projection onto $\textrm{span}\{e^{in\f}\}_{n\ge 0}$. The APS boundary conditions on~$D$  amount to the following choice of the domain:
\begin{gather}\label{ddomain}
\textrm{dom}(D) = \left\{b\in\mathcal{H} \, : \  \|Db\|^2_w<\infty,\   P_{\ge 0}\sigma(b)=0\right\}.
\end{gather}
There are certain subtleties in this def\/inition which are clarif\/ied in the statement of Proposition~\ref{ddecomp} at the end of this section.

The next proposition describes a (partial) Fourier series decomposition of the Hilbert space~$\mathcal{H}$. Def\/ine
\begin{gather*}
a(k) := \frac{w(k)^2}{S(k)},
\end{gather*}
and let
\begin{gather*}
\wlltwoz = \left\{\{g(k)\}_ {k\in\Z}\, : \ \|g\|^2_a=\sum_{k\in\Z}a(k)^{-1}|g(k)|^2 <\infty\right\}.
\end{gather*}
Now we are ready for the Fourier decomposition of $\mathcal{H}$ which is just like~(\ref{hdecomp1}).

\begin{proposition}
Let $\mathcal{H}$ be the Hilbert space defined above.  Then the formula
\begin{gather*}
b = \sum_{n\in\Z} g_n(K)\left(U^*\right)^n
\end{gather*}
defines an isomorphism of Hilbert spaces
\begin{gather}\label{hdecomp}
\bigoplus_{n\in\Z} \wlltwoz \cong \mathcal{H}.
\end{gather}
\end{proposition}

\begin{proof}
The proof is identical to the one in \cite{KM}. In particular we have
\begin{gather*}
\|b\|^2_w = \sum_{n\in\Z} \textrm{tr}\left(S(K)w^{-2}(K)|g_n(K)|^2\right).\tag*{\qed}
\end{gather*}
\renewcommand{\qed}{}
\end{proof}

The main reason we consider the Fourier decomposition is that it (again partially) diagonalizes the operator $D$. This is the subject of the next lemma. Before stating it we need some more notation. Consider the ratios:
\begin{gather*}
c^{(n)}(k) := \frac{w(k+n)}{w(k)}
\end{gather*}
and notice that since $\{w(k)\}$ is an increasing sequence we have:
\begin{alignat*}{3}
& c^{(n)}(k)=1 \qquad&& \textrm{for} \ n=0, & \\
& c^{(n)}(k)>1 \qquad&& \textrm{for} \ n>0, & \\
& c^{(n)}(k)<1 \qquad&& \textrm{for} \ n<0. &
\end{alignat*}
The coef\/f\/icients $c^{(n)}$ are needed to def\/ine the following operators in~$\wlltwoz$. The f\/irst is:
\begin{gather*}
\overline{A}^{(n)}g(k) = a(k)\big(g(k) - c^{(n)}(k)g(k+1)\big)
\end{gather*}
with domain
\begin{gather*}
\textrm{dom}(\overline{A}) = \left\{g\in\wlltwoz \, : \ \|\overline{A}g\|_{a}<\infty \right\}.
\end{gather*}
Additionally consider the operator
$\overline{A_0}^{(n)}$ which is the operator $\overline{A}^{(n)}$ but with domain
\begin{gather*}
\textrm{dom}\big(\overline{A_0}^{(n)}\big) = \big\{g\in \textrm{dom}\big(\overline{A}^{(n)}\big)\, : \ g_\infty:=\lim_{k\to\infty}g(k) = 0 \big\} .
\end{gather*}
The last def\/inition makes sense since by the analysis of~\cite{KM} the limit $\lim_{k\to\infty}g(k)$ exists for $g\in \textrm{dom}(\overline{A})$.
One has the following proposition, which is a quantum analog of the f\/irst part of Theorem~\ref{theo1}.

\begin{proposition}\label{ddecomp}
With respect to the decomposition \eqref{hdecomp} one has:
\begin{gather*}
D \cong \bigoplus_{n>0} \overline{A}^{(n)} \oplus \bigoplus_{n\leq 0} \overline{A_0}^{(n)}.
\end{gather*}
Equivalently:
\begin{gather*}
Db = \sum_{n>0} \overline{A}^{(n)}g_n(K)(U^*)^n+\sum_{n\leq 0} \overline{A_0}^{(n)}g_n(K)(U^*)^n,
\end{gather*}
where
\begin{gather*}
b = \sum_{n\in\Z} g_n(K)\left(U^*\right)^n
\end{gather*}
\end{proposition}

\begin{proof}
The proof is a direct calculation identical to the one in~\cite{KM}.
\end{proof}

\section{Construction  of the parametrix}\label{section4}
In this section we construct and analyze in detail the inverse (= a parametrix) $Q$ for the  ope\-ra\-tor~$D$.   The construction is fairly similar to the one done in Section~4 in~\cite{KM}, however the norm estimates are quite dif\/ferent. Somewhat surprisingly the norm estimates below hold for any choice of sequence of weights $\{w(k)\}$ satisfying~(\ref{weightconditions}).

We start with a lemma containing estimates of sums through integrals. Recall that the  sequence~$\{w(k)\}$ is increasing with limits at $\pm\infty$ equal, correspondingly, to~$w^+$ and~0.

\begin{lemma}\label{intineq}
If $f(t)$ is a decreasing continuous function on $(0,(w^+)^2)$ then
\begin{gather}\label{ineq1}
\sum_{l<k}f\big(w(k)^2\big)S(k) = \sum_{l<k}f\big(w(k)^2\big)\big(w(k)^2-w(k-1)^2\big) \le \int_{w(l)^2}^{w_+^2} f(t)dt,\\
\label{ineq2}
\sum_{k\le l}f\big(w(k)^2\big)S(k) \le \int_0^{w(l)^2} f(t)dt,\\
\label{ineq3}
\sum_{k\in\Z}f\big(w(k-1)^2\big)S(k) \ge \int_0^{w_+^2} f(t)dt.
\end{gather}
\end{lemma}

The proof of the statements of the lemma follows from a straightforward comparison of the Riemann sums of the left hand side with the integrals on the right hand side.

Our presentation in this section is as follows. First we discuss the kernels of the $\overline{A}^{(n)}$ operators for the three cases $n=0$, $n>0$, $n<0$.   Secondly we construct the parametrices for all three cases.   Thirdly we discuss the norm estimates of the parametrices, and f\/inally we summarize the analysis in the main result of this paper.

Below we show that the operator $D$ has no kernel by analyzing the terms in the decomposition of Proposition~\ref{ddecomp}.

\begin{proposition}\label{kernelnzero}
The operators  $\overline{A}^{(n)}$ for $n\ge 0$ and $\overline{A_0}^{(n)}$ for $n<0$ have no kernel.
\end{proposition}

\begin{proof}
We start with $n=0$. Here $c^{(n)}(k)=1$ and
it is clear that the only solution of that $\overline{A}^{(0)}R^{(0)} = 0$ is, up to a constant, $R^{(0)} = 1$.   But one has
\begin{gather*}
\big\|R^{(0)}\big\|^2_a  = \sum_{k\in\Z}\frac{1}{a(k)}=\sum_{k\in\Z}\frac{S(k)}{w(k)^2}\ge \textrm{const}\sum_{k\in\Z}\frac{S(k)}{w(k-1)^2}=\infty,
\end{gather*}
where we used condition~4 of~(\ref{weightconditions}) as well as~(\ref{ineq3}) for $f(t)=1/t$.
Therefore $R^{(0)}(K) \not\in \wlltwoz$ and hence $\overline{A}^{(0)}$ has no kernel.

Next we discuss the kernel of $\overline{A}^{(n)}$ when $n>0$. It is not too hard to see that any element of the kernel has to be proportional to
\begin{gather*}
R^{(n)}(k) := \prod_{l=k}^\infty c^{(n)}(l) = \frac{(w_+)^n}{w(k)w(k+1)\cdots w(k+n-1)}.
\end{gather*}
The norm calculation gives
\begin{gather*}
\|R^{(n)}\|^2_a = \sum_{k\in\Z}\frac{1}{a(k)} |R^{(n)}(k)|^2 = \sum_{k\in\Z}\frac{1}{a(k)} \prod_{l=k}^\infty |c^{(n)}(l)|^2.
\end{gather*}
Since $|c^{(n)}(l)|>1$ and $\sum_{k\in\Z}\frac{1}{a(k)}=\infty$, the  sum above diverges and hence $R^{(n)}(K) \not\in\wlltwoz$.

Finally we discuss the kernel of the operator $\overline{A_0}^{(n)}$ when $n<0$.   Yet again the kernel is formally one dimensional
and spanned by
\begin{gather*}
R^{(n)}(k) = \prod_{l=k}^\infty c^{(n)}(l) = \frac{w(k+n)w(k+n-1)\cdots w(k-1)}{(w_+)^{-n}}.
\end{gather*}
While one can easily show that $R^{(n)} \in \wlltwoz$,  one however has $\lim\limits_{k\to\infty} R^{(n)}(k) = 1 \neq 0$, so this means $R^{(n)} \not\in \textrm{dom}\big(\overline{A_0}^{(n)}\big)$.   Thus the result follows.
\end{proof}

The second portion of the discussion  is the construction of the parametrices for all three cases.  Since there are no kernels (and cokernels) involved
we simply compute the inverses of opera\-tors~$\overline{A}^{(n)}$. Thus, given $g(k)$, one needs to solve the equation $\overline{A}^{(n)}f(k) = g(k)$
where additionally we need $\lim\limits_{k\to\infty}f(k) = 0$ for $n\le 0$. This is done in a similar manner to the methods in \cite{CKW,KM}.   In the case when $n>0$ one arrives at the following formula:
\begin{gather*}
f(k) = -\sum_{l<k} \frac{R^{(n)}(k)}{R^{(n)}(l)a(l)}g(l) = -\sum_{l<k}\frac{w(l)\cdots w(l+n-1)}{w(k)\cdots w(k+n-1)} \frac{S(l)}{w(l)^2}g(l).
\end{gather*}
Similarly in the case $n\le 0$ one has:
\begin{gather*}
f(k) = \sum_{k\le l}\frac{R^{(n)}(k)}{R^{(n)}(l)}  \frac{g(l)}{a(l)} = \sum_{k\le l}\frac{w(k+n)\cdots w(k-1)}{w(l+n)\cdots w(l-1)} \frac{S(l)}{w(l)^2}g(l).
\end{gather*}

The right hand sides of the above equation give the parametrices $Q^{(n)}$ for all three cases.   We thus have the following:
\begin{alignat}{3}
&Q^{(n)}g(k) = -\sum_{l<k}\frac{S(l)}{w(l)^2}g(l) \qquad && \textrm{for} \ n=0, & \nonumber\\
&Q^{(n)}g(k) = -\sum_{l<k}\frac{w(l)\cdots w(l+n-1)}{w(k)\cdots w(k+n-1)} \frac{S(l)}{w(l)^2}g(l)\qquad && \textrm{for} \ n>0, & \label{qnformulas}\\
&Q^{(n)}g(k) = \sum_{k\le l}\frac{w(k+n)\cdots w(k-1)}{w(l+n)\cdots w(l-1)}  \frac{S(l)}{w(l)^2}g(l)\qquad && \textrm{for} \ n<0. & \nonumber
\end{alignat}

We summarize the above analysis in the following proposition.

\begin{proposition}\label{qn}
Let $Q^{(n)}$ be defined by the formulas above, then we have the following
\begin{alignat*}{4}
&\overline{A}^{(n)}Q^{(n)} = I \qquad && \textrm{and}\quad Q^{(n)}\overline{A}^{(n)} = I \qquad && \textrm{for }n>0; & \\
&\overline{A_0}^{(n)}Q^{(n)} = I\qquad && \textrm{and}\quad Q^{(n)}\overline{A_0}^{(n)} = I \qquad && \textrm{for }n\le 0. &
\end{alignat*}
\end{proposition}

Next we discuss the boundedness for the parametrices in the cases $n>0$ and $n<0$.
The dif\/f\/iculty comes for $k\to\infty$: while the ratios of weights are always less than 1, the series $\sum_{k\in\Z}\frac{S(k)}{w(k)^2}$ is not summable and we cannot replicate the estimates of~\cite{CKW} and~\cite{KM}. In fact the integral operators $Q^{(n)}$ are not Hilbert--Schmidt. The trick is to estimate most but not all weight ratios by one. The remaining sums, containing potentially divergent terms, are estimated by integrals using Lemma~\ref{intineq}.
We have the following result.

\begin{proposition}\label{boundednessofQ}
The operators $Q^{(n)}$ defined above are bounded operators in~$\wlltwoz$ when $n\ne 0$.
\end{proposition}

\begin{proof}
First consider the case that $n>0$.  Applying the Schur--Young inequality and the inequalities~(\ref{ineq1}),
and (\ref{ineq2}) one has
\begin{gather*}
\big\|Q^{(n)}\big\|_a^2
\le \sup_{k} \left(\sum_{l<k}\frac{w(l)\cdots w(l+n-1)}{w(k)\cdots w(k+n-1)} \frac{S(l)}{w(l)^2}\right)\sup_{l} \left(\sum_{l<k}\frac{w(l)\cdots w(l+n-1)}{w(k)\cdots w(k+n-1)} \frac{S(k)}{w(k)^2}\right) \\
\phantom{\big\|Q^{(n)}\big\|_a^2 }{}
 \le \sup_{k} \left(\frac{1}{w(k)}\sum_{l<k}\frac{S(l)}{w(l)}\right)\sup_{l} \left(w(l)\sum_{l<k}\frac{S(k)}{w(k)^3}\right) \\
\phantom{\big\|Q^{(n)}\big\|_a^2 }{}
 \le \sup_{k} \left(\frac{1}{w(k-1)}\sum_{l\le k-1}\frac{S(l)}{w(l)}\right)\sup_{l} \left(w(l) \int_{w(l)^2}^{w_+^2} t^{-\frac{3}{2}}dt\right) \\
\phantom{\big\|Q^{(n)}\big\|_a^2 }{}
 \le \sup_{k} \left(\frac{1}{w(k-1)} \int_0^{w(k-1)^2}t^{-\frac{1}{2}}dt\right)\cdot2 \sup_{l} \left(1-\frac{w(l)}{w_+}\right) \le 2\cdot 2 = 4.
\end{gather*}

Thus $Q^{(n)}$ is bounded for $n>0$.   Next consider the case $n<0$.   Here one has quite similar estimates:
\begin{gather*}
 \big\|Q^{(n)}\big\|_a^2  \le \sup_{k} \left(\sum_{k\le l}\frac{w(k+n)\cdots w(k-1)}{w(l+n)\cdots w(l-1)} \frac{S(l)}{w(l)^2}\right)\sup_{l} \left(\sum_{k\le l}\frac{w(k+n)\cdots w(k-1)}{w(l+n)\cdots w(l-1)} \frac{S(k)}{w(k)^2}\right) \\
 \phantom{\big\|Q^{(n)}\big\|_a^2}{}
 \le \sup_{k} \left(w(k-1)\sum_{k\le l}\frac{S(l)}{w(l)^2w(l-1)}\right)\sup_{l} \left(\frac{1}{w(l-1)}\sum_{k\le l}\frac{S(k)}{w(k)}\right) \\
\phantom{\big\|Q^{(n)}\big\|_a^2}{}
\le \left(\sup_{l} \frac{w(l)}{w(l-1)}\right)\sup_{k} \left(w(k-1)\sum_{k-1 < l}\frac{S(l)}{w(l)^3}\right)\sup_{l} \left(\frac{1}{w(l-1)}\sum_{k\le l}\frac{S(k)}{w(k)}\right) \\
\phantom{\big\|Q^{(n)}\big\|_a^2}{}
\le \left(\sup_{l} \frac{w(l)}{w(l-1)}\right)\sup_{k} \left(w(k-1) \int_{w(k-1)^2}^{w_+^2}t^{-\frac{3}{2}}dt\right)\sup_{l} \left(\frac{1}{w(l-1)}\int_0^{w(l)^2}t^{-\frac{1}{2}}dt\right) \\
\phantom{\big\|Q^{(n)}\big\|_a^2}{}
\le 4\left(\sup_{l} \frac{w(l)}{w(l-1)}\right)^2 <\infty.
\end{gather*}
Thus $Q^{(n)}$ is bounded for $n<0$ and this completes the proof.
\end{proof}

Finally we put together the previous information about the parametrix $Q$ of the Dirac ope\-ra\-tor~$D$ def\/ined in Section~\ref{section3}.   We state the main result of this note.

\begin{theorem}
Let $D$ be the operator \eqref{ddef} with domain \eqref{ddomain}.  Then there exists an ope\-ra\-tor~$Q$ such that $QD=DQ=I$. Moreover, with respect to the decomposition \eqref{hdecomp}
one has
\begin{gather}\label{qdecomp}
Q = \bigoplus_{n\in\Z} Q^{(n)}=Q^{(0)}+\tilde{Q},
\end{gather}
where the operators $Q^{(n)}$ are given by \eqref{qnformulas} and $\tilde{Q}$ is bounded.
\end{theorem}

\begin{proof}
By Proposition~\ref{ddecomp} one can decompose $D$ as $\bigoplus_{n>0} \overline{A}^{(n)} \oplus \bigoplus_{n\leq 0} \overline{A_0}^{(n)}$
which in turn gives the decomposition~(\ref{qdecomp}) of $Q$.
One has that
\begin{gather*}
\|\tilde Q\|_{w} = \sup_{n\ne 0} \big\|Q^{(n)}\big\|_a.
\end{gather*}
Then from Proposition~\ref{boundednessofQ}, one has the following inequalities
\begin{gather*}
\|Q\|^2_{w} \le 4\left(\sup_{l} \frac{w(l)}{w(l-1)}\right)^2 < \infty,
\end{gather*}
where the last inequality follows from the assumptions in~(\ref{weightconditions}).   To see that one has $DQ{=}QD{=}I$ we use the decompositions of~$Q$ and $D$ and Proposition~\ref{qn}. This completes the proof.
\end{proof}

\section{The balanced quantum Dirac operators}\label{section5}
In this section we study a version of the constructions of the previous sections that is more like the theory of \cite{KM}.
The main objects: the Hilbert space and the Dirac operator are called balanced since in their def\/initions the left multiplication is not preferred over the right multiplication.

Since the results for the balanced Dirac operators are completely analogous to the ``un\-ba\-lan\-ced'' case and the proofs require only trivial modif\/ication, we simply state the main steps of the construction. The only signif\/icant dif\/ference between the two cases are the estimates on the components of the parametrix.

In order to avoid unnecessary complications we simply recycle the old notation. As before the starting point is the choice of a sequence of weights $\{w(k)\}_{k\in\Z}$ satisfying~(\ref{weightconditions}).   The Hilbert space $\mathcal{H}$ is the space of power series:
\begin{gather*}
b = \sum_{n\in\Z} g_n(K)\left(U^*\right)^n
\end{gather*}
but this time with a dif\/ferent, balanced norm:
\begin{gather*}
\|b\|^2_w  =  \textrm{tr}\left(S^{1/2}w(K)^{-1}bb^*w(K)^{-1}S^{1/2}\right) \\
\phantom{\|b\|^2_w }{} =\sum_{n\in\Z} \textrm{tr}\left(\sqrt{S(K)S(K+n)}w^{-1}(K)w^{-1}(K+n)|g_n(K)|^2\right).
\end{gather*}

The balanced Dirac operator is
\begin{gather*}
Db = -S^{-1/2}U^*\left[b,U_w\right]w(K)S^{-1/2}
\end{gather*}
with the domain:
\begin{gather*}
\textrm{dom}(D) = \left\{b\in\mathcal{H} \, : \ \|Db\|^2_w<\infty,\  P_{\ge 0}\sigma_{\textrm{circ}}(b)=0\right\}.
\end{gather*}

As before the Dirac operator splits into Fourier components. To describe them we need to modify the coef\/f\/icients of the previous sections.
Actually, the coef\/f\/icients
\begin{gather*}
c^{(n)}(k) := \frac{w(k+n)}{w(k)}
\end{gather*}
stay the same, but we need to change:
\begin{gather*}
a^{(n)}(k) := \frac{w(k)w(k+n)}{\sqrt{S(k)S(k+n)}}.
\end{gather*}
Those are used for the following previously def\/ined operators in $\wlltwoz$. The f\/irst operator is:
\begin{gather*}
\overline{A}^{(n)}g(k) = a^{(n)}(g(k) - c^{(n)}(k)g(k+1))
\end{gather*}
with domain
\begin{gather*}
\textrm{dom}(\overline{A}) = \left\{g\in\wlltwoz \, : \ \|\overline{A}g\|_{\wlltwoz}<\infty \right\},
\end{gather*}
and the second operator
$\overline{A_0}^{(n)}$ is the operator $\overline{A}^{(n)}$ but with domain
\begin{gather*}
\textrm{dom}\big(\overline{A_0}^{(n)}\big) = \big\{g\in \textrm{dom}\big(\overline{A}^{(n)}\big) \, : \  g_\infty:=\lim_{k\to\infty}g(k) = 0 \big\} .
\end{gather*}

With that notation, the Proposition~\ref{ddecomp} remains true. In particular one has:
\begin{gather*}
D \cong \bigoplus_{n>0} \overline{A}^{(n)} \oplus \bigoplus_{n\leq 0} \overline{A_0}^{(n)}.
\end{gather*}

The problem of inverting the operator $D$ is tackled as in the previous section. The components of the inverse are given by formulas like
(\ref{qnformulas}) with the only modif\/ication coming from the dif\/fe\-rent~$a^{(n)}$ coef\/f\/icients. We end up with the following expressions for the parametrices:
\begin{alignat*}{3}
&Q^{(n)}g(k) = -\sum_{l<k}\frac{\sqrt{S(l)S(l+n)}}{w(l)w(l+n)}g(l)\qquad && \textrm{for} \ n=0, & \\
&Q^{(n)}g(k) = -\sum_{l<k}\frac{w(l)\cdots w(l+n-1)}{w(k)\cdots w(k+n-1)} \frac{\sqrt{S(l)S(l+n)}}{w(l)w(l+n)}g(l)\qquad && \textrm{for}\ n>0, & \\
&Q^{(n)}g(k) = \sum_{k\le l}\frac{w(k+n)\cdots w(k-1)}{w(l+n)\cdots w(l-1)}  \frac{\sqrt{S(l)S(l+n)}}{w(l)w(l+n)}g(l)\qquad&& \textrm{for} \ n<0. &
\end{alignat*}
One can verify directly that for the operator $Q = \bigoplus_{n\in\Z} Q^{(n)}$ we have $QD=DQ=I$.
The following is the main result of this section.

\begin{proposition}
The operators $Q^{(n)}$ defined above are bounded operators in $\wlltwoz$ when $n\ne 0$.
\end{proposition}

\begin{proof}
We use the Schur--Young inequality and follows the steps of the proof of the Proposi\-tion~\ref{boundednessofQ}, with some modif\/ications.
We show the details for $n<0$, the other case is completely analogous.

There are two sums that we need to estimate. The f\/irst sum is:
\begin{gather*}
\Sigma^n_1(k):=\sum_{k\le l}\frac{w(k+n)\cdots w(k-1)}{w(l+n)\cdots w(l-1)} \frac{\sqrt{S(l)S(l+n)}}{w(l)w(l+n)}.
\end{gather*}
Using Cachy--Schwartz inequality we estimate:
\begin{gather*} \Sigma^n_1(k)\le
 \left(\sum_{k\le l}\frac{w(k+n)\cdots w(k-1)}{w(l+n)\cdots w(l-1)} \frac{S(l)}{w(l)^2}\right)^{1/2}\!\!
\left(\sum_{k\le l}\frac{w(k+n)\cdots w(k-1)}{w(l+n)\cdots w(l-1)} \frac{S(l+n)}{w(l+n)^2}\right)^{1/2} \\[-0.5ex]
\phantom{\Sigma^n_1(k)}{}  \le \left(w(k-1)\sum_{k\le l}\frac{S(l)}{w(l)^2w(l-1)}\right)^{1/2}\left(w(k+n)\sum_{k\le l}\frac{S(l+n)}{w(l+n)^3}\right)^{1/2}.
\end{gather*}
Since the weights in the denominator are bigger than the corresponding weights in the numerator, their ratios were estimated by one. The f\/irst term on the rights hand side of the above was already estimated
in the proof of Proposition~\ref{boundednessofQ}. The second term is essentially the same as the f\/irst:
\begin{gather*}
\sup_{k} \left( w(k+n)\sum_{k\le l}\frac{S(l+n)}{w(l+n)^3}\right)=
\sup_{k} \left(w(k)\sum_{k\le l}\frac{S(l)}{w(l)^3}\right).
\end{gather*}
It follows that $\Sigma^n_1(k)$ is bounded uniformly in $n$.

The second sum in the Schur--Young inequality is
\begin{gather*}
\Sigma^n_2(l):=\sum_{k\le l}\frac{w(k+n)\cdots w(k-1)}{w(l+n)\cdots w(l-1)} \frac{\sqrt{S(k)S(k+n)}}{w(k)w(k+n)}
\end{gather*}
and we bound it in the same fashion as the f\/irst sum:
\begin{gather*}
\Sigma^n_2(l)
 \le\left(\sum_{k\le l}\frac{w(k+n)\cdots w(k-1)}{w(l+n)\cdots w(l-1)} \frac{S(k)}{w(k)^2}\right)^{1/2}
\left(\sum_{k\le l}\frac{w(k+n)\cdots w(k-1)}{w(l+n)\cdots w(l-1)} \frac{S(k+n)}{w(k+n)^2}\right)^{1/2} \\[-0.5ex]
 \phantom{\Sigma^n_2(l) }{}
 \le \left(\frac{1}{w(l-1)}\sum_{k\le l}\frac{S(k)}{w(k)}\right)^{1/2}\left(\frac{1}{w(l+n)}\sum_{k\le l}\frac{S(k+n)}{w(k+n)}\right)^{1/2}.
\end{gather*}
Again the f\/irst term above was already estimated
in the proof of Proposition~\ref{boundednessofQ}, and the second term is essentially the same as the f\/irst.
It follows that $\Sigma^n_1(k)$ is uniformly bounded. Repeating the same steps for $n>0$ gives the boundedness of $Q$
for the balanced Dirac operator.
\end{proof}

\pdfbookmark[1]{References}{ref}
\LastPageEnding

\end{document}